\documentclass[12pt,twoside]{amsart}
\usepackage{amsmath}
\usepackage{amssymb}
\usepackage{amscd}
\usepackage{xypic}
\xyoption{all}
\title[A note on Asymptotic Serre duality for jet differentials]{A note on Asymptotic Serre duality for jet differentials}

\author{Mohammad Reza Rahmati}
\thanks{}
\address{ Abdus Salam School of Mathematical Sciences , GCU, Lahore, PAKISTAN 
\hfill\break 
\hfill\break \\
\hfill\break }
\email{rahmati@sms.edu.pk}

\newcommand{\comments}[1]{}

\newtheorem{theorem}{Theorem}[section]

\keywords{Green-Griffiths jet bundle, Serre duality, Entire Curve Locus, Holomorphic 
Morse Inequalities, Weighted Projective Space, Relative Canonical Sheaf}

\subjclass{32Q45, 30D35, 14E99}

\begin{document}

\maketitle

\begin{abstract}
We formulate and prove the existence of an asymptotic duality along the fibers of the Green-Griffiths jet bundles over 
projective manifolds. The existence of global sections for these bundles and also for their dual sheaves has been established by 
the author and J. P. Demailly in a previous article, \cite{DR}.
\end{abstract}



\section{Introduction}

The main result of this paper is a Serre duality for asymptotic sections of jet bundles as formulated in Theorem 2.1 below. An 
application is also given for partial resolution of the Green-Griffiths conjecture.

P. Griffiths and M. Green \cite{GG} provided several equivalent methods for characterizing the bundle of jets over a complex
projective 
variety $X$. Later, the techniques developed by Semple, Demailly, McQuillan, Siu and Green-Griffiths led to new conjectures in 
K\"ahler geometry. Demailly, \cite{D1}, \cite{D2}, shows the existence of global sections of the Green-Griffiths bundles with 
sufficiently high degree.

Any entire curve $ f: \mathbb{C} \to X $ satisfies $Q(f_{[k]})=0$, where $f_{[k]}$ a lift of $f$ and $Q$ a global section of the 
jet bundle (see \cite{D1}). As the section $Q$ involves both the manifold and jet coordinates, it is desirable to show that 
there are enough dual global sections along the jet fibers to pair with the given section to give equations on the variety $X$.

The Green-Griffiths bundle associated with a pair $ (X, V) $, where $ X $ is a projective variety and $ V \subset T_X $ is a 
holomorphic sub-bundle of rank $r$, is a sequence of projective bundles over $X$ which are defined inductively as in
\cite[page 20]{D1}. Alternatively, one can define the Green-Griffiths bundle $X_k$ as
$$
X_k: = (J_kV \smallsetminus {0}) / \mathbb{C}^* \, ,
$$
where $J_k$ is the bundle of germs of $k$-jets of Taylor expansion for $f$. The bundles $\pi_k:X_k \to X$ provide a tool to study entire holomorphic curves $ f: \mathbb{C} \to X $, since such a curve has a lift to $f_{[k]}:\mathbb{C} \to  X_k $ for every $ k $. Using the notation in \cite{D1}, we write $E_{k, m}V^* =(\pi_{k})_* \mathcal{O}_{X_k}(m)$. 

In \cite{D1}, \cite{D2}, the existence of asymptotic global sections for Green-Griffiths jets on a complex projective variety $X$ has 
been proved, by using Morse estimates for the curvature of suitable metrics on these bundles. It is shown in \cite{D1} that
$$
H^0(X,E_{k,m}V^* \otimes A^{-1})=H^0(X_k,\mathcal{O}_{X_k}(m) \otimes \pi_{k}^*A^{-1}))
$$
are non-trivial when $m \gg k \gg 0$, where $A$ is a positive line bundle on $X$. In fact the factor $A$ can
be absorbed in the other over $X_k$ if we assume $m \in \mathbb{Q}$. 

The formulation of Serre duality is based on the existence of the canonical sheaf. A technical difficulty arises here that in 
the definition of the relative canonical sheaf, the fibers of $X_k \to X$ have singularities. However, there is a generalization 
of the definition of the canonical sheaf for singular varieties as explained in \cite{D1}, \cite{DR}.

The relative canonical sheaf $K_{X_k/X}$ can be defined in a manner similar to the definition of the canonical sheaf $K_V$, 
where $V$ is a holomorphic subbundle of $T_X$. Serre duality can then be written as
$$
H^0(X_k,(\pi_k)_*\mathcal{O}_{X_k}(m)) \bigotimes H^{k(r-1)}(X_k, K_{X_k/X} \otimes (\pi_k)_*\mathcal{O}_{X_k}(-m)) \longrightarrow \mathcal{O}_X
$$
for $ m \gg k \gg 0 $. The non-triviality of the cohomology groups has been proved by Demailly and the author, \cite{D1}, 
\cite{D2} and \cite{DR}.

\section{Serre duality for Jet bundles}

As explained above the existence of the relative canonical sheaves $K_{X_k/X}$ provides the possibility to establish Serre 
duality along fibers of jet bundles. Because the fibers in the Green-Griffiths bundle are weighted projective spaces of 
appropriate weights $(1,2,\cdots ,n)$, the relative Serre duality can be interpreted as the duality for coherent sheaves on 
weighted projective spaces. On account of this we review the formulation of the adjoint pair along the fibers of $E_{k,m}V^*$. 
Then, according to the classical Serre duality the dual pair associated to $H^{0} (F_x,(\pi_k)_* \mathcal{O}_{{X_k},x} (m))$ is 
$H^{k (r-1)} (F_x, K_{F_x} \otimes (\pi_k)_*\mathcal{O}_{{X_k},x} (- m))$, where $ K_{F_x} $ is the canonical sheaf of the 
fiber. In other words we have
$$
H^0(\pi_{k}^{-1}(x),(\pi_k)_* 
\mathcal{O}_{X_k}(m))^{\vee}=H^{k(r-1)}(\pi_{k}^{-1}(x), K_{F_x} \otimes (\pi_k)_*\mathcal{O}_{X_k}(-m)).
$$

By the Leray spectral sequence for $\pi_k:X_k \to X$ we get the following:
$$
H^{k (r-1)} (X_k, K_{X_k / X} \otimes (\pi_k)_*\mathcal{O}_{X_k} (- m ))
$$
$$
=\,
H^0 (X, R^{k (r-1)} (\pi_{k})_* (K_{X_k / X} \otimes (\pi_k)_*\mathcal{O}_{X_k} (-m))).
$$

\noindent
On account of the formalism of Serre duality for coherent sheaves on projective spaces, the candidate for the duality is
$$
H^0(X_k,(\pi_k)_*\mathcal{O}_{X_k}(m)) \bigotimes H^{k(r-1)}(X_k, K_{X_k/X} \otimes (\pi_k)_*\mathcal{O}_{X_k}(-m)) \longrightarrow \mathcal{O}_X\, .
$$
We have to ensure that the sheaf $H^{k(r-1)}(X_k, K_{X_k / X} \otimes (\pi_k)_*\mathcal{O}_{X_k} (- m)) $ is 
non-trivial i.e. has enough sections for $ m \gg k \gg 0 $.

\begin{theorem}[{Serre Duality for Jet Fibers}]
There is a Serre duality for asymptotic cohomologies, $m \gg k \gg 0$,
\begin{equation}
H^0(X_k,(\pi_k)_*\mathcal{O}_{X_k}(m)) \bigotimes H^{k(r-1)}(X_k, K_{X_k/X} \otimes (\pi_k)_*\mathcal{O}_{X_k}(-m))
\longrightarrow \mathcal{O}_X.
\end{equation} 
\end{theorem}

\begin{proof} The duality in construction is the duality on each fiber of the jet bundle, all glued by the spectral sequence of Leray. In that sense, it is a duality of coherent sheaves on the weighted projective spaces \cite{RT}. The adjoint pair on the fiber along the fiber $\pi_k^{-1}(x)$ is given by the formula 
$$
H^0(\pi_{k}^{-1}(x),(\pi_k)_* \mathcal{O}_{X_k}(m))^{\vee}=
H^{k(r-1)}((\pi_{k}^{-1}(x), K_{F_x} \otimes (\pi_k)_*\mathcal{O}_{X_k}(-m)).
$$
The degeneration of the Leray spectral sequence of the fibration $\pi_k:X_k \to X$ provides 
$$
H^0(X_k,(\pi_k)_*\mathcal{O}_{X_k}(m))^{\vee} = H^0(X, R^{k(r-1)}(\pi_{k})_*( K_{X_k/X} \otimes (\mathcal{O}_{X_k}(-m)))\, , 
$$
which is equivalent to (1). The non-triviality of the first factor in (1) is proved in \cite{D1}. The nontriviality of the adjoint cohomology group in the pairing is obtained from estimates
$$
H^q(X_k, K_{X_k/X} \otimes (\pi_k)_*\mathcal{O}_{X_k}(-m))\geq \sum_{q-1,q,q+1}\frac{rm^n}{r!}\int_{X(\Theta,j)}(-1)^{q-j}\Theta^n-o(m^n)\, ,
$$
where $\Theta$ is the curvature of suitable k-jet metric on $X_k$ and $$X(\Theta,j)=\{x \in X; \ \Theta \ 
\text{has signature} \ (n-j,j)\}\, .$$ We have considered the inequality for $q=k(r-1)$. It follows that
when the $ K_{X_k / X} $ is big, both of the factors in the pairing (1) are nontrivial for $m\gg 0$, cf.
\cite[Section 4]{DR}. The theorem follows.
\end{proof}

In \cite{M}, J. Merker shows that when $X$ is a hypersurface of degree $d$ in $ \mathbb{P}^{n + 1}$ and is a generic member of 
the universal family $\mathfrak{X} \subset \mathbb{P}^{n+1} \times \mathbb{P}^{N_d}$, the Green-Griffiths conjecture holds for 
$X$. His method uses ideas of Y. T. Siu, on existence of slanting vector fields, see \cite{S}, \cite{P}. The proof establishes 
the conjecture outside a certain algebraic subvariety $ \Sigma \subset J_{\text{vert}}^n (\mathfrak{X}) $ defined by Wronskians. 
An implication of Merker's result \cite{M} is the following:
$$
H^0 \big (\mathfrak{X}_k,(\pi_k)_*\mathcal{O}_{\mathfrak{X}_k}(m) \big ) \bigotimes H^{k(r-1)}\big (\mathfrak{X}_k, K_{\mathfrak{X}_k/\mathfrak{X}} \otimes \big \langle J_{\text{vert}}^{k} (\mathfrak{X}) \big \rangle^{m} \big ) \longrightarrow \mathcal{O}_\mathfrak{X} \, ,
$$
where $ m \gg k \gg 0$, and $\langle J_{\text{vert}}^{k} (\mathfrak{X}) \rangle^{m}$ means the ring of operators generated 
by the $J_{\text{vert}}^{k} (\mathfrak{X})$ of degree $m$.

\vspace{0.5cm}

\noindent
\textit{Application to Green-Griffiths Conjecture:} We note that any entire curve $f:\mathbb{C} \to X$ is satisfied by the
sections $Q \in H^0(X_k, (\pi_k)_*\mathcal{O}_{X_k} (m)) $ i.e., $Q(f_{[k]})=0$ for $f_{[k]}:\mathbb{C} \to X_k$ a lift of
$f$, cf. \cite{D1}. It follows that, in (1) after the composition for $Q \in H^0 (X_k,(\pi_k)_*\mathcal{O}(X_k (m) )$ and
$P \in H^0 (X_k,K_{X_k/X} \otimes (\pi_k)_*\mathcal{O}(X_k (-m) )$, the resulting function also satisfies
$$
\langle Q,P \rangle (f_{[k]})=0, \qquad Q=\sum_{|\alpha|=m}A_{\alpha}(z)\xi^{\alpha}, \ P= \sum_{|\beta|=m'}B_{\beta}(z)\partial_{\xi}^{\beta}.
$$

Therefore, the issue is to show that the image of (1) is a non-trivial ideal of $\mathcal{O}_X$. The above procedure plays a 
similar task as the idea on the existence of slanting vector fields along jet fibers, due to Siu \cite{S}. In fact the slanting 
vector fields play the role of generating sections in the adjoint fiber rings of jet bundles.

\vspace{0.5cm}

\noindent
{ACKNOWLEDGEMENT:} I thank the Abdus Salam School of Mathematical Sciences , GCU, for its research facilities and financial support.

\noindent
{CONFLICT OF INTEREST:} There is no conflict of interest.

 \end{document}